\documentclass[12pt]{amsart}

\usepackage{amsgen,amsmath,amstext,amsbsy,amsopn,amsfonts,amssymb}

\def\Proof{\noindent{\sl Proof.}\ }
\def\qed{{\hfill $\Box$ \medbreak}}

\newtheorem{defi}{Definition}[section]
\newtheorem{lem}[defi]{Lemma}

\newtheorem{cor}[defi]{Corollary}

\newtheorem{prop}[defi]{Proposition}
\newtheorem{example}[defi]{Example}

\newtheorem*{thmm}{Main Theorem}

\DeclareMathOperator{\la}{\langle}
\DeclareMathOperator{\ra}{\rangle}

\DeclareMathOperator{\ad}{ad}

\DeclareMathOperator{\F}{\mathbb{F}}

\begin{document}

\title[Lie metabelian restricted enveloping algebras]{Restricted enveloping algebras whose skew and symmetric elements are Lie metabelian}
\author{\textsc{Salvatore Siciliano}}
\address{Dipartimento di Matematica e Fisica ``Ennio De Giorgi", Universit\`{a} del Salento,
Via Provinciale Lecce--Arnesano, 73100--Lecce, Italy}
\email{salvatore.siciliano@unisalento.it}

\author{\textsc{Hamid Usefi}}
\address{Department of Mathematics and Statistics,
Memorial University of Newfoundland,
St. John's, NL,
Canada, 
A1C 5S7}
\email{usefi@mun.ca}

\thanks{The research of second author was supported by NSERC of Canada. }

\begin{abstract}Let $L$ be a restricted Lie algebra over a field of characteristic $p>2$ and denote by $u(L)$ its restricted enveloping algebra.
We establish when the symmetric or skew  elements of $u(L)$ under the principal involution are Lie metabelian.
\end{abstract}
\subjclass[2010]{16S30,  17B60, 16R40, 16W10,  17B50}
\date{\today}


\keywords {Restricted Lie algebra, enveloping algebra, polynomial identity, Lie metabelian, skew and symmetric elements}

\maketitle

\section{Introduction}
Let $A$ be an algebra with involution $\ast$ over a field $\F$. We denote by $A^+=\{x\in A\vert\,x^\ast=x\}$ the set of symmetric elements of $A$ under $\ast$ and by $A^-=\{x\in A\vert\,x^\ast=-x\}$ the set of skew-symmetric elements.  A general question of interest is to establish the extent to which the structure of $A^+$ or $A^-$ determines the structure of $A$ (see \cite{H76}).  For instance, a celebrated result of Amitsur in \cite{A68} states that if $A^+$ or $A^-$ satisfies a polynomial identity, then so does $A$. Moreover, a considerable amount of attention has been devoted to decide if Lie properties satisfied by  the symmetric or the skew symmetric elements of a group algebra $\F G$ under the canonical involution are also satisfied by the whole algebra $\F G$, see e.g. 
\cite{GPS1, GS93, JR, LSS, LSS1}.

Now, let $L$ be a restricted Lie algebra  over a field ${\mathbb F}$ of characteristic $p>2$ and let $u(L)$ be the restricted enveloping algebra of $L$.  
We denote by $\top$  the \emph{principal involution} of $u(L)$, that is, the unique ${\mathbb F}$-antiautomorphism of $u(L)$ such that $x^\top=-x$ for every $x$ in $L$. We recall that $\top$ is just the antipode of the  ${\mathbb F}$-Hopf algebra $u(L)$.
In \cite{Skew} and \cite{SUIsr} the conditions under which $u(L)^-$ or $u(L)^+$ are Lie solvable, Lie nilpotent or bounded Lie Engel were provided. It turns out that 
$u(L)^-$ or $u(L)^+$ are Lie solvable if and only if so is $u(L)$.
The aim of this note is to characterize $L$ when $u(L)^-$ or $u(L)^+$ are Lie metabelian.  Our main result is the following:

\begin{thmm}\label{skew-metabelian}
Let $L$ be a restricted Lie algebra over a field $\F$ of characteristic $p>2$. 
Then the following hold:
\begin{enumerate}
\item [{\normalfont 1)}] $u(L)^-$  is Lie  metabelian if and only if  either $L$   is   abelian or $p=3$, $L^\prime$ is 1-dimensional and central, and ${L^\prime}^{[p]}=0$.
\item [{\normalfont 2)}] $u(L)^+$  is Lie  metabelian if and only if  one of the following conditions is satisfied:
 {\begin{enumerate}
\item [{\normalfont (i)}] $L$ is abelian;
\item [{\normalfont (ii)}] $p=3$, $L^\prime$ is 1-dimensional and central, and ${L^\prime}^{[p]}=0$;
\item [{\normalfont (iii)}] $p=3$ and $L$ is 2-dimensional. 
\end{enumerate}}
\end{enumerate}
\end{thmm}

Note that Lie metabelian restricted enveloping algebras have been characterized in  \cite{S3}. By combining this result and our main theorem, one concludes that in odd characteristic $u(L)^-$ is Lie metabelian if and only if so is $u(L)$. This remains true for the symmetric case provided that $p>3$, but if $L$ is a 2-dimensional non-abelian restricted Lie algebra over a field of characteristic $3$, then $u(L)^+$ is Lie metabelian whereas $u(L)$ is not. It seems interesting that this is  indeed the only exception. We also show that in characteristic 2 our main theorem  fails both for  skew and symmetric case. We finally mention that analogous results for group algebras have been carried out in 
\cite{BRR,CLS,LS, LR}.



\section{Proof of the Main Theorem}\label{} 
Throughout the paper, unless otherwise stated, $\F$ will  denote a field of characteristic $p\geq 3$.

Let $A$ be an associative algebra over $\F$. 
The Lie bracket on $A$ is defined by $[x,y]=xy-yx$, for every $x,y\in A$.
We set $[x_1,x_2]^o =[x_1,x_2]$ and define inductively
$$
[x_1,x_2,\ldots,x_{2^{n+1}}]^o=[[x_1,\ldots,x_{2^n}]^o,
[x_{2^n+1},\ldots,x_{2^{n+1}}]^o].
$$ 
A subset $S$ of $A$ is said to be Lie solvable if there exists an $n$ such that 
$$
[x_1,x_2,\ldots,x_{2^{n+1}}]^o=0,
$$
 for every $x_1,\ldots,x_{2^{n+1}}\in S$. In particular, if $n=2$ then we say $S$ is Lie metabelian.

If $L$ is a restricted Lie algebra over  $\F$ then we denote by $L'$ the derived subalgebra of $L$.
For a subset $S$ of $L$ we denote by $\la S\ra_p$
the restricted ideal of $L$ generated by $S$. The centralizer of $S$ in $L$ is denoted by $C_L(S)$ and $Z(L)$ is the center of $L$.  Moreover, for every positive integer $n$, we use  $\zeta_n(L)$ to denote 
the $n$-th term of the ascending central series of $L$.

For the proof of  our main result we first need to prove some technical lemmas.

\smallskip
\begin{lem}\label{nonheis} 
 Let $L$ be a restricted Lie algebra and suppose that $L$ contains an $\F$-linearly independent set
  $\{ a,b,c,v,w\}$ such that
$[a,b]=v$, $[a,c]=w$, $[b,c]=0$, and $v,w\in Z(L)$.
 Then $u(L)^-$ and $u(L)^+$ are not Lie  metabelian.
\end{lem}
\Proof Since all the elements $a,b,a^2w,c^2v$ are skew-symmetric, we have
$$
[[a^2w,b],[a,c^2v]]=[2avw,2cvw]=4vw^2\neq 0,
$$
by the PBW Theorem for restricted Lie algebras (see e.g. \cite[\S 2, Theorem 5.1]{SF}). Hence, $u(L)^-$ is not Lie metabelian. 
Now note that $2ac-w$ and $2ab-v$ are symmetric and so we have
$$
[[2ac-w, c^2], [av, 2ab-v]]=[ 4c^2 w, 2av^2]=-16 cv^2w^2\neq 0,
$$
by the PBW Theorem. Hence,  $u(L)^+$ is not Lie metabelian. \qed

\begin{lem}\label{heis} 
Let $L$ be a restricted Lie algebra  and suppose that $L$ contains an $\F$-linearly independent set
 $\{ x_1,x_2,y_1,y_2,v,w\}$ such that
$[x_1,x_2]=v$, $[y_1,y_2]=w$, $[x_i,y_j]=0$ for every $i,j\in \{1,2\}$, and $v,w\in Z(L)$.
Then $u(L)^-$ and $u(L)^+$ are not Lie  metabelian.
\end{lem}
\Proof
Since $x_1,x_2,x_1y_1w,x_2y_2^2 \in u(L)^-$, we have 
$$
[[x_1y_1w,x_2],[x_1,x_2y^2]]=[y_1vw,y_2^2v]=2y_2v^2w^2\neq 0,
$$ 
by the PBW Theorem. Hence,  $u(L)^-$ is not Lie metabelian. As for the symmetric case, note that
$$
[ [x_1y_1, x_2 y_1], [x_1y_1, y_2^2]]=[ y_1^2v, 2x_1y_2w]=4x_1y_1v w^2\neq 0,
$$
implying that $u(L)^+$ is not Lie metabelian. \qed

\begin{lem}\label{class3} 
 Let $L$ be a restricted Lie algebra and suppose that $L$ contains an $\F$-linearly independent set
 $\{ x,y,v,w\}$ such that
$[x,y]=v$, $[y,v]=w$, $[x,v]=0$, and $w\in Z(L)$.
Then $u(L)^-$ and $u(L)^+$ are not Lie  metabelian. 
\end{lem}
\Proof  Note that $u=2xyw-vw$ is skew-symmetric. Since
$$
[[u,y],[x,y]]=[2vyw+w^2,v]=2vw^2\neq 0,
$$
we deduce that $u(L)^-$ is not Lie metabelian.
Since all the elements $x^2,y^2,vw,2xy-v$ are symmetric, we have
$$
[[x^2,2xy-v],[y^2,vw]]=[4x^2v,2yw^2]=16xv^2w^2-8x^2w^3\neq 0,
$$
by the PBW Theorem, hence $u(L)^+$ is not Lie metabelian. \qed

The following elementary result is likely well-known. However, since we do not have a reference, we offer a short proof. 
\begin{prop}\label{locally nilpotent}
Let $L$ be locally nilpotent Lie algebra over any field. If $L^\prime$ is finite-dimensional then $L$ is nilpotent of class at most $\dim L^\prime+1$.
\end{prop}
\Proof It is enough to show that every finite-dimensional subalgebra $H$ of $L$ is nilpotent of class at most 
$\dim L^\prime+1$. 
To do so, we prove  by induction on $n=\dim H^\prime$ that $H$ is nilpotent of class at most $n+1$.
If $n=0$, the assertion is clear. Now let $H$ be a non-abelian finite-dimensional subalgebra of $L$. 
Since $H$ is nilpotent, there exists a non-zero central element $x\in H^\prime$. Now consider 
$\bar H=H/\la x\ra_{\F}$.  Since $\dim \bar H^\prime<n$, we deduce by the induction hypothesis that
$\bar H$ is nilpotent of class at most $n$. Hence, $H$ is nilpotent of class at most $n+1$, as required.
\qed

\begin{lem}\label{nilpotent}
Let $L$ be a restricted Lie algebra. If $u(L)^-$  is Lie  metabelian then  $L$ is nilpotent. The same conclusion holds when $u(L)^+$ is Lie metabelian provided $p\neq 3$.
\end{lem}
\Proof We know by Theorem 1 in \cite{Skew} and Theorem 1.3 in \cite{SUIsr} that $L^\prime$ is finite-dimensional.  Hence, by Proposition \ref{locally nilpotent}, it    is enough to show that $L$ is locally nilpotent.  Without loss of generality, we can assume that the ground field is algebraically closed. Now, let $H$ be a finite-dimensional subalgebra of $L$ and assume that $H$ is not nilpotent. By the Engel's Theorem, then there exists an element $y\in H$ such that the adjoint map
$\ad y$ is not a nilpotent linear transformation. Hence, $\ad y$ has a non-zero eigenvalue $\lambda$. Thus,  there exists $x\in H$ such that $[x,y]=\lambda x$. Now we rescale $y$ to assume that $[x,y]= x$. 
Note that the element $2xy^2-2xy+x$ is skew-symmetric and one has 
$$
[[x,y], [2xy^2-2xy+x,y]]=[x, 2xy^2-2xy+x]=4x^2y-4x^2\neq 0,
$$
by the PBW Theorem. Therefore, if $u(L)^-$   is Lie  metabelian, we have a contradiction. Hence, $H$ must be nilpotent. 

For the symmetric case note that $2xy-x$ is symmetric and we have
\begin{align*}
[[ 2xy-x, y^2], [ x^2, 2xy-x]]&=[ 4xy^2-4xy+x, 4x^3]\\
&=48(-x^3y+3x^3+x^4)\neq 0,
\end{align*}
noting that $p>3$ and  by the PBW Theorem. 
 Hence, $u(L)^+$   is not Lie  metabelian which is a contradiction.
We conclude again that $H$ must be nilpotent.
\qed

In characteristic 3, $L$ need not be nilpotent when $u(L)^+$ is Lie metabelian. Indeed, we have:

\begin{lem}\label{2-dimensional}
Let $L$ be the 2-dimensional non-abelian restricted Lie algebra over a field $\F$ of characteristic  $3$. Then $u(L)^+$  is Lie  metabelian.
\end{lem}
\Proof It is easy to see that $L$ has a basis  $x,y$ such that  $[x,y]=x$, $x^{[3]}=0$, $y^{[3]}=y$ (see Section 2.1 of \cite{SF}).  As $\F$ has odd characteristic, $u(L)^+$ is spanned by the trace elements $a+a^\top$, $a\in u(L)$. Thus we have
$$
u(L)^+=\langle 1, 2xy-x, x^2,y^2, x^2y^2-2x^2y\rangle_{\F}.
$$ 
By explicit computations we get 
$$
[u(L)^+,u(L)^+]=\langle x^2y-x^2, xy^2-xy+x\rangle_{\F}.
$$
As $[x^2y-x^2, xy^2-xy+x]=0$, we conclude that $u(L)^+$  is Lie  metabelian.  \qed

\begin{lem}\label{abelian}
Let $L$ be a restricted Lie algebra. If $u(L)^+$  or $u(L)^-$  is Lie  metabelian then  either  $L$  is   abelian or $p=3$ and  the power mapping acts trivially on central commutators of $L$.
\end{lem}
\Proof Let $x,y\in L$ and set $z=[x,y]$. Suppose that $z$ is central. We have
\begin{align*}
[[2xy-z, yz], [ x^2,2xy-z]]&= [2yz^2, 4x^2z ]=-16xz^4\in \delta_2(u(L)^+);\\
[[x^2z, y], [ x^2y-xz, y]]&= [2xz^2, 2xyz-z^2 ]=4xz^4\in \delta_2(u(L)^-).
\end{align*}
Hence, if $u(L)^+$  or $u(L)^-$ is Lie  metabelian then we deduce by the PBW Theorem that either $z=0$ or $p=3$ and $z^3=0$. We conclude that  if $u(L)^+$  or $u(L)^-$  is Lie  metabelian then either $L$ is abelian or $p=3$ and the power mapping acts trivially on central commutators.
\qed


\begin{lem}\label{nonnilp}
Let $L$ be a non-nilpotent restricted Lie algebra over a field $\F$ of characteristic $3$. If $u(L)^+$  is Lie  metabelian then $L$  is 2-dimensional.
\end{lem}
\Proof  Without loss of generality, we can assume that the ground field $\F$ is algebraically closed.   In view of Theorem 1.3 of \cite{SUIsr}, $L^\prime$ is finite-dimensional. 
Since $L$ is not nilpotent we deduce from  Proposition \ref{locally nilpotent}  that $L$ contains a non-nilpotent finite-dimensional subalgebra $H$.  By the Engel's Theorem, there exists $y\in H$ such that the the adjoint map $\ad y$ is not nilpotent. Let $\lambda$ be a non-zero eigenvalue of $\ad y$. Then there exists $x\in H$ such that $[x,y]=\lambda x$. Now replace $y$ by $\lambda^{-1}y$ to assume that $[x,y]= x$. 

We claim that $C_L(x,y)=0$. Let $a\in C_L(x,y)$. As $u(L)^+$ is Lie metabelian we must have
$$
0=[[xa,ya],[2xy-x,y^2]=[xa^2,xy^2-xy+x]=2x^2a^2y-2x^2a^2,
$$
and so the PBW Theorem forces that $a=\alpha x+\beta y$ for some $\alpha, \beta \in \F$. As $[a,x]=[a,y]=0$, it follows that $a=0$, as claimed.
 
Since $\dim L/C_L(x)=\dim [L,x]$ and $\dim L/C_L(y)=\dim [L,y]$  are both finite, we conclude that $L$ is finite-dimensional. Moreover, as $x^{[3]}$ and $y-y^{[3]}$ both commute with $x$ and $y$, this also entails that $x^{[3]}=0$ and $y=y^{[3]}$. Therefore we have that $\ad y=(\ad y)^3$ and so $L$ decomposes as 
$$
L=L_0\oplus L_1 \oplus L_2,
$$
 where $L_\lambda$ denotes the eigenspace corresponding to the eigenvalue $\lambda$ of the linear transformation $\ad y$, for $\lambda=0,1,2$.
 
We now claim that $L_0=\F y$. Suppose to the contrary and let $a\in L_0$ such that $a$ and $y$ are linearly independent. If we had that $[a,x]=k x$ for some $\lambda \in \F$, then $[x,a+k y]=0=[y,a+k y]$ and so $a+k y \in C_L(x,y)=0$, a contradiction.  As $[a,x]\in L_1$, it follows that the set $\{x,y,a,[a,x]\}$ is $\F$-linearly independent. Put $b=[[a,x],x]\in L_2$ and note that $[b,x]=[a,x^{[3]}]=0$. As $u(L)^+$ is Lie metabelian, we have
$$
0=[[x^2,y^2],[y^2,b^2]]=[x^2y-x^2,b^2y+b^2]=x^2b^2y.
$$ 
Hence, $x,y$ and $b$ are linearly dependent by the PBW Theorem. Since $x,y$ and $b$ are in distinct eigenspaces of $\ad x$, we conclude that $b=0$. Put $c=[a,x]\in L_1$ and suppose, if possible, that $c\neq 0$. As $x$ and $c$ are linearly independent, by the PBW Theorem, we have
$$
[[xc,y^2],[2xy-x,2cy-c]]=[xcy-xc,2xc]=2x^2c^2 \neq 0,
$$ 
hence $u(L)^+$ is not metabelian, a contradiction. Consequently, we have $a\in C_L(x,y)=0$, another contradiction, which yields the claim.

Next we prove that $L_2=0$. Suppose by way of contradiction that $L_2$ contains a non-zero element $z$. Then $[z,x]\in L_0$ and so, for what was proved above, we have $[x,z]=\beta y$, for some $\beta\in \F$. 
As $u(L)^+$ is Lie metabelian, we have
\begin{align*}
0&=[[x^2,y^2],[y^2,2zy+z]]=[x^2y-x^2,zy^2+zy+z]\\ &=2x^2zy^2-x^2zy+(1-\beta)xy^2+(1-\beta)xy.
\end{align*}
As the elements $x,y,$ and $z$ are linearly independent, this contradicts the PBW Theorem, yielding the claim.

We finally prove that $L_1=\F x$. Suppose to the contrary and let $v\in L_1$ such that $v$ and $x$ are linearly independent. Then $[v,x]\in L_2=0$ and  by the PBW Theorem we have
$$
[[x^2,y^2],[y^2,v^2]]=[x^2y-x^2,2v^2y+v^2]=2x^2v^2\neq 0,
$$
contradicting the fact that $u(L)^+$ is Lie metabelian. This finishes the proof. \qed

We are now in a position to prove our main result:

\medskip
\emph{Proof of the Main Theorem.} 
The sufficiency of both parts of the statement follows from \cite{S3}. Let us prove the necessity.
Note that, by Lemma \ref{abelian}, if $u(L)^+$  is  Lie metabelian  and $L$ is not nilpotent then  case 2(iii) occurs. Therefore, for the rest of the proof, we assume that 
$u(L)^+$ or $u(L)^-$ is  Lie metabelian and $L$ is nilpotent. 
Furthermore, by Lemma \ref{abelian},  it is enough to show that if $p=3$ and $L$ is not abelian then $L'$ is 1-dimensional.  Assume, by contradiction, that 
$\dim_{\F}  L^\prime>1$. 
Let $n$ be the nilpotence class of $L$ and put ${\mathfrak L}=L/\zeta_{n-2}(L)$. We  proceed by considering the following three  cases.
  
\emph{ Case 1:} $\dim_{\F} {\mathfrak L}^\prime>1$ and there exist $x_1,x_2,x_3\in {\mathfrak L}$ such that the elements $x_4=[x_1,x_2]$ and $x_5=[x_1,x_3]$ are $\F$-linearly independent.
 
If $x_6=[x_2,x_3]\notin  \la x_4, x_5\ra_{\F}$ then since
${\mathfrak L}$ is nilpotent of class 2, we have that $\{x_1,x_2,x_3,x_4,x_5,x_6 \}$ is $\F$-linearly independent. 
Furthermore, from Lemma \ref{abelian} applied to ${\mathfrak L}$  we have that 
 $\la x_6\ra_{\F}$ is a restricted ideal of $L$. Put 
$$
\bar{H}=\langle x_1,x_2,x_3,x_4,x_5,x_6 \rangle_p/\la x_6\ra_{\F}.
$$  
 Now we apply Lemma \ref{nonheis} to $\bar{H}$ and get a contradiction. 

On the other hand, if $x_6=\alpha x_4 + \beta x_5$ for some $\alpha, \beta \in \F$, then put 
$$
\tilde{H}=\langle x_1,x_2,x_3,x_4,x_5 \rangle_p 
$$
and, for every $i=1,2,\dots,5$, define:
$$
\tilde{x_i}=
\begin{cases}
-\beta x_1+x_2, &\text{if $i=2$};\\
\alpha x_1+x_3, &\text{if $i=3$};\\
x_i, &\text{otherwise}.
\end{cases}
$$
Clearly, $\tilde{H}$ and $\{\tilde{x_1},\tilde{x_2},\tilde{x_3},\tilde{x_4},\tilde{x_5} \}$ satisfy the hypotheses of Lemma \ref{nonheis} and we get a contradiction again.
 
\emph{ Case 2:} $\dim_{\F} {\mathfrak L}^\prime>1$ and for every $a,b,c\in {\mathfrak L}$ the elements $[a,b]$ and $[a,c]$ are $\F$-linearly dependent. Let 
$x_1,x_2,y_1,y_2 \in {\mathfrak L}$ such that the elements 
$v=[x_1,x_2]$ and $w=[y_1,y_2]$ are $\F$-linearly independent and 
 observe that $x_i$ and $y_j$ must commute, for  $i,j =1,2$. Since ${\mathfrak L}$ is nilpotent of class 2, 
 $\{x_1,x_2,y_1,y_2,v,w\}$ is $\F$-linearly independent. 
 Now we apply Lemma \ref{heis} to 
$
 {H}=\langle x_1,x_2,y_1,y_2,v,w \rangle_p
$
 and get a  contradiction.

\emph{ Case 3:} $\dim_{\F}  {\mathfrak L}^\prime=1$. Since $\dim_{\F}  L^\prime\geq 2$,
$L$ has nilpotency class at least 3. Put ${\bar L}=L/\zeta_{n-3}(L)$. Note that
$$
\dim_{\F} \frac{{\bar L}^\prime+Z({\bar L})}{Z({\bar L})}=1.
$$
It follows that $\bar{L}$ is  not 2-Engel because every 2-Engel  Lie algebra is nilpotent of class at most 2.
As a consequence, we can find elements $x,y\in \bar{L}$ such that 
$
[[x,y],y]\neq 0.
$
Put $v=[x,y]$, $w=[y,v]$, and $z=[x,v]$.  
Suppose first $z\notin \la w\ra_{\F}$. Note that  $\{x,y,v,w,z \}$ is $\F$-linearly
independent. Moreover, by Lemma \ref{abelian}, $\la z\ra_{\F}$ is a restricted ideal. Now we apply 
Lemma \ref{class3} to 
$
\bar {H}=\langle x,y,v,w,z \rangle_p/\la z\ra_{\F}$   to get a contradiction.

Finally if $z=\alpha w$ for some $\alpha \in \F$ 
then replace $x$ by $x-\alpha y$
and apply Lemma \ref{class3} to
$
H=\langle x,y,v,w \rangle_p$ to get a contradiction again, which completes the proof.
\qed 

By combining Theorem \ref{skew-metabelian} with the main result of \cite{S3} we get:

\begin{cor}\label{cor-metabelian}
Let $L$ be a restricted Lie algebra over a field of characteristic $p>2$.  Then $u(L)^-$  is Lie  metabelian if and only if  so is $u(L)$.
\end{cor}

Unlike the  skew case, in  view of Lemma \ref{nonnilp},   the fact that the symmetric elements are Lie metabelian does not force  the whole algebra $u(L)$ is  Lie metabelian if $p=3$.   
We also observe that Corollary \ref{cor-metabelian} and our main theorem fail  in characteristic 2, as the following example shows.

\begin{example}
\emph{
Let $\F$ be a field of characteristic 2 and consider the restricted Lie algebra $L$ over $\F$ given by  $L=\la x, y, z\mid [x,y]=x, [x,z]=[y,z]=z^2=0, x^2=z, y^2=y\rangle$. 
Note that the $p$-map is not trivial on $L'$. However,  it is easy to see
that $u(L)^+=u(L)^-=\la x, y, z, xy, yz\ra_{\F}$. Hence, $u(L)^+=u(L)^-$ is Lie metabelian despite the fact that 
$u(L)$ is not.}
\end{example}

\enddocument